\magnification 1200

\input amstex
\documentstyle{amsppt}

\def\ho#1#2{\hbox{$H_{#1}(#2)$}}
\def\bd{\partial}
\def\cp{\hbox{${\Bbb C}P^2$}}
\def\ncp{\hbox{${\overline {{\Bbb C}P}^2}$}}
\def\r{\hbox{$\Bbb R$}}
\def\rf{\hbox{${\Bbb R}^4$}}
\def\br#1{\hbox{${\Cal R}_{#1}$}}

\topmatter
\title{Smooth structures on collarable ends of 4-manifolds} \endtitle

\author{\v Zarko Bi\v zaca and John Etnyre} \endauthor

\address{Department of Mathematics, University of Texas, Austin, TX 78712}
\endaddress
\email{bizaca\@math.utexas.edu, etnyre\@math.utexas.edu} \endemail

\date April 18, 1996\enddate

%%\subjclass \endsubjclass
%%\keywords \endkeywords 

\abstract We use Furuta's result, usually referred to as ``10/8-conjecture'',
to show that for any compact 3-manifold $M$ the open manifold $M\times\r$ has
infinitely many different smooth structures.  Another consequence of Furuta's
result is existence of infinitely many smooth structures on open topological
4-manifolds with a topologically collarable end, provided there are only
finitely many ends homeomorphic to it. We also show that for each closed spin
4-manifold there are exotic \rf's that can not be smoothly embedded into
it.\endabstract
\endtopmatter

\document
\baselineskip = 13pt
\parskip = 4pt

There is a well known (and still unsolved at the time of writing) conjecture in
four dimensional topology known as the {\it 11/8 conjecture}.  Namely, if
$X$ is a closed spin topological four-manifold, then for an appropriately
chosen basis for its second homology the intersection matrix has the form
$(\oplus_kE_8)\oplus (\oplus_lH)$, where $E_8$ is the negative definite 8 by 8
intersection matrix of the ``$E_8$ plumbing'', and $H$ is a 2 by 2 hyperbolic
matrix,

$$H=\pmatrix 0 & 1\cr 1 & 0\cr\endpmatrix.$$

\noindent
We consider the case when $X$ is a {\it smooth} manifold.  By Rochlin's theorem
the number of copies of $E_8$, $k$, is even.  Further, the 11/8 conjecture
asserts that the number of hyperbolics, $l$, in the intersection form of $X$,
satisfies the relation

$$l\geq \frac{3}{2}k.$$

\noindent
The numbers in the ``11/8'' refer to the quotient of the Euler characteristic
and the signature of the manifold in the in the case of the equality
$l=\frac{3}{2}k$.  The conjecture is known to be true in the case $k=2$, see
Theorem 2 of [D], where the minimal possible $l$ is realized by the K3
surface.  Recently Furuta has proven the ``10/8'' inequality, that is, for a
given number of $E_8$'s, $k$, the number of hyperbolics, $l$, has to satisfy
the inequality $l > k+1$, see [S].

In this work we use Furuta's result, or more precisely a weaker inequality
$l\geq k$, to prove the existence of infinitely many different smooth
structures on $M\times\r$, where $M$ is an arbitrary compact three-manifold.
Furthermore, we show that 4-manifolds with a topologically collarable end that
satisfy the conditions of Theorem 5 have infinitely many smooth structures.
Our construction of these smooth structures is based on Gompf's [G1] end-sums
with exotic \rf's.  The main building block in our constructions is a well
known exotic \rf, which we will denote by \br1, and can be obtained from the
K3 surface.  Our smooth structures are obtained as end-sums of a given open
manifold with different number of copies of \br1 and, consequently, we are
producing only countable families of smooth structures.  The infinite end-sum
of copies of \br1 and their mirror images is reminiscent to the universal
exotic \rf\ constructed by Freedman and Taylor [FT], although our construct is
probably not ``universal'' itself, it has the same known non embedding
properties of the universal \rf.

The authors would like to thank R.~Gompf and D.~Freed for their helpful
suggestions.

\subhead{} Construction and properties of \br{n}\endsubhead 
Let $Y$ denote the standard K3 surface.  Its intersection form introduces a
decomposition of the second homology, $\ho{2}{Y}\cong (\oplus_2E_8)\oplus
(\oplus_3H).$ According to Casson [C] the six elements of \ho{2}{Y} that span
the $\oplus_3H$ summand can be represented by three wedges of two immersed
2-spheres and the extra intersections can be grouped into pairs consisting of
one with a positive and one with a negative sign.  For each of these
$\pm$-pairs of intersections there is a loop with a framed neighborhood such
that if a 2-handle is attached to it with the prescribed framing then there is
an ambient isotopy inside the union of a neighborhood of the 2-spheres and the
2-handle that eliminates the given pair of intersections.  Although, in
general, such a 2-handle can not be embedded ambiently, Casson has shown that
we can ambiently cap the loop by a Casson handle.  We fix an open regular
neighborhood of the immersed wedges of 2-spheres and denote its union with the
capping Casson handles by $N$. Since any Casson handle can be embedded into the
standard 2-handle, $N$ smoothly embeds into $\sharp_3(S^2\times S^2)$.  From
Freedman's celebrated result [F] it follows that each Casson handle has a core
that is locally flat topologically embedded 2-disc. The complement of the
wedges of 2-spheres and the (topological) cores of the Casson handles in
$\sharp_3(S^2\times S^2)$ is an exotic \rf, that is, a manifold homeomorphic
but not diffeomorphic to the standard \rf.  We denote this manifold by \br1.
We select an arbitrary 3-manifold inside $N$ that separates the end of $N$ from
the wedges of spheres and topological core discs. It also separates
$\sharp_3(S^2\times S^2)$ and the closure of the component contained in \br1 we
denote by $K_1$. So, $K_1$ is a compact manifold with boundary contained in
\br1. Note that there are neighborhoods of ends of \br1 and the interior of
$K_1$ that are smoothly embedded into $Y$ and that each has a separating
topologically embedded locally flat 3-sphere.

For an integer $n>1$ we define \br{n} to be end-sum of $n$ copies of \br1, $\
\br{n}=\natural_n\br1$.  Each copy of \br1 inside \br{n} contains a copy of
$K_1$.  We use the connecting neighborhoods of arcs from the end-sum to form a
boundary connected sum of the copies of $K_1$, $K_n = \natural_n K_1$, so
$K_n$ is a compact manifold with boundary embedded into \br{n}.  Equivalently,
we can construct \br{n} by working with a connected sum of $n$ copies of the
K3 surface: let $Y_n$ denote connected sum of $n$ copies of the K3 surface, $\
\ho{2}{Y_n}\cong(\oplus_{2n}E_8)\oplus (\oplus_{3n}H)$. Each copy of the K3
surface in $Y_n$ contains a copy of $N$. We assume that the 4-balls affected by
the construction of the connected sums were disjoint from the copies of $N$, so
we connect together the copies of $N$ by open regular neighborhood of arcs,
each passing through one of the $S^3\times I$ tubes that connect the K3
surfaces in $Y_n$.  We call the resulting open connected manifold $N_n$.  As in
the case of single K3 surface, we can smoothly transplant $N_n$ into
$\sharp_{3n}(S^2\times S^2)$. It is easy to see that the complement of the
wedges of 2-spheres and topological cores in $\sharp_{3n}(S^2\times S^2)$ is
diffeomorphic to \br{n}. Again, as in the case of single K3 surface, there is a
neighborhood $A_n$ of the end of the interior of $K_n$ that is contained in
$N_n$ and, therefore, in $Y_n$.  We note that there is a topological locally
flat embedding of the 3-sphere that separates $A_n$.

We also define $\br{\infty}=\natural_{\infty}\br1$ and $\
\br{*}=\br{\infty}\natural-(\br{\infty})$, where $-(\br{\infty})$ denotes
the ``mirror image'' of \br{\infty}, the same underlying manifold, but with
opposite orientation.

\proclaim{Theorem 1} If $X$ is a closed smooth spin 4-manifold, then there is
an integer $m>0$ such that for any integer $n\geq m$ \br{n} and $K_n$ can not
be smoothly embedded into $X$ or into any closed smooth spin 4-manifold with
the same intersection form as $X$.
\endproclaim

\demo{Proof} As we have noted above, $\ \ho{2}{X}\cong(\oplus_{2k}E_8)\oplus 
(\oplus_lH)$.  Let $m$ be any integer with the property $2m > l-2k$.  We
claim that when $n\geq m$ \ \br{n} cannot be smoothly embedded into $X$.  Let
us assume that there is such an embedding.  From the image of this embedding we
remove the complement of $A_n$ in $K_n$ and denote the resulting open manifold
by $X_0$. The end of $X_0$ has a neighborhood diffeomorphic to $A_n$.  We have
constructed above an embedding of $A_n$ via $N_n$ into $Y_n$, the connected sum
of $n$ copies of the K3 surface.  We cut out from $Y_n$ the component of $N_n -
A_n$ that contains the wedges of immersed 2-spheres and denote the resulting
manifold $Y_{n,0}$.  We form a closed smooth manifold $Z$ by identifying the
ends of $X_0$ and $Y_{n,0}$.  Since the left out piece contained the generators
of the second homology summand, $\ho{2}{Y_{n,0}}\cong\oplus_{2k}E_8$. Note
that using the topological locally flat 3-sphere separating $A_n$ it
is easy to see that $\ho{2}{Z}\cong\ho{2}{X}\oplus\ho{2}{Y_{n,0}}$, so
$$\ho{2}{Z}\cong(\oplus_{2k+2n}E_8)\oplus (\oplus_lH)$$
with $2k+2n> l$ which is impossible by Furuta's result.
\qed\enddemo

An obvious consequence of this theorem is 

\proclaim{Corollary 2} For $1\leq m,n\leq \infty$ and $n\neq m$, \br{n} 
and \br{m} are not diffeomorphic. \br{\infty} and \br{*} can not be smoothly
embedded into any closed smooth spin 4-manifold.  \br{n}, $K_n$, where $n>0$,
\br{\infty} and \br{*} can not be embedded into any negative definite smooth
4-manifold.  Furthermore, $\br{n}\natural -(\br{n})$ and \br{*} do not embed
into any definite closed smooth 4-manifold.
\endproclaim

Note that since the universal \rf\ from [FT] contains any exotic \rf, it itself
can not be embedded into a closed smooth spin 4-manifold. 

\demo{Proof} From the construction it follows that \br{n} can be embedded into 
$\sharp_{3n}(S^2\times S^2)$, but not into $\sharp_{2n}(S^2\times S^2)$;
otherwise, as in the proof of Theorem 1 we could construct a closed spin
4-manifold that would contradict Furuta's result.  We take \br{n_1} to be
\br{1}.  Let us assume that we have chosen nondiffeomorphic \br{n_i}, for
$i\leq k$.  We choose an integer $n_{k+1} > \frac{3}{2}n_k$.  Then \br{n_{k+1}}
does not smoothly embed into $\sharp_{3n_{k}}(S^2\times S^2)$ and so it is not
diffeomorphic to any \br{n_i}, for $i\leq k$. Therefore, there are infinitely
many smoothly distinct \br{n}'s. Moreover, no two of \br{n}'s are
diffeomorphic, since each \br{n} is end sum of $n$ copies of \br{1}, the
equation $\br{n}=\br{m}, \ n<m,$ would imply that for any $k> m$,
$\br{k}=\br{l}$ for some $n\leq l\leq m$.  So it would follow that there are
only finitely many distinct \br{n}'s.

Since \br{\infty} and \br{*} contain every \br{n} it follows from Theorem 1
that \br{\infty} does not embed into any spin closed smooth 4-manifold. Thus it
follows that \br{\infty} can not be diffeomorphic to any \br{n}, for $1\leq n
<\infty$.

If \br{n} or $K_n$ could be embedded into a negative definite smooth
4-manifold, we could excise all hyperbolics from the connected sum of $n$ K3
surfaces and reglue the complement of the embedded \br{n} inside a negative
definite manifold.  The resulting closed 4-manifold would be a smooth negative
definite 4-manifold whose intersection form would contain $E_8$ summands.  It
is well known algebraic result that such forms can not be diagonalized over the
ring of integers so the existence of such a manifold is impossible by
Donaldson's theorem (see Theorem 1 in [D], to remove the restriction on the
fundamental group one needs to use Seiberg-Witten theory).  Therefore, \br{n},
for $n>0$, and so \br{\infty} and \br{*} can not be embedded into any
negative definite simply-connected smooth 4-manifold.

Finally, since $-\br{n}$ does not embed into a positive definite smooth
4-manifold, $\br{n}\natural -(\br{n})$ and \br{*} do not embed into any
definite closed smooth 4-manifold.\qed\enddemo

\remark{Remark} Following an argument of Freedman, see [G1], we may find in
each \br{n} uncountable many nondiffeomorphic exotic \rf's sharing the same
nonembedding property as \br{n}.  Namely, in \br{n} there is a neighborhood of
the end that is homeomorphic to $S^3\times (0,\infty)$ and disjoint from
$K_n$. For each $0 <t<\infty$ let \br{n,t} denotes the complement of $S^3\times
[t,\infty)$.  It follows from the work of Taubes [T] on periodic ends that no
two of \br{n,t}, $0 <t<\infty$, are diffeomorphic.
\endremark

$\br{*}$ shares another non embedding property with the universal \rf\ from
[FT]: it can not be embedded into the interior of a smooth 4-cell.

\proclaim{Proposition 3} If $e: B^4\hookrightarrow X$ is a locally flat 
topological embedding of the 4-ball into a smooth 4-manifold $X$, then \br{*}
can not be smoothly embedded into the interior of $e(B^4)$.
\endproclaim

\demo{Proof} We follow Freedman and Taylor's proof of the same property for the
universal exotic \rf. The first step is to embed $X$ into $\
\cp\sharp_n(S^2\times S^2)$, for some $n$ large enough.  So if $X$ had a
boundary we double it, then surger homotopicaly non trivial loops away from
$e(B^4)$ and form the connected sum with sufficient numbers of copies of \cp\
and \ncp\ to obtain a closed smooth simply connected manifold with the required
intersection form.  By adding copies of $S^2\times S^2$ the manifold becomes
diffeomorphic to $\cp\sharp_n(S^2\times S^2)$.  The second homology class of
a complex projective line inside $\cp\subset\cp\sharp_n(S^2\times S^2)$ can
be represented by an embedded surface in the complement of $e(B^4)$.  By adding
extra copies of $S^2\times S^2$ if necessary, we can find in the complement of
$e(B^4)$ a smooth 2-sphere representing the same homology class.  Blowing down
this sphere produces an embedding of $e(B^4)$ into $\sharp_n(S^2\times S^2)$.
It follows from Corollary 2 that there can not be a smooth embedding of \br{*}
into $e(B^4)$.\qed\enddemo

\subhead{} Smooth structures on $M^3\times\r$ \endsubhead 
Gompf has used his end-sum construction to produce uncountably many smooth
structures on $S^3\times\r$.  Recently F. Ding [Di] has extended this result
to produce an uncountable family of smooth structures $M^3\times\r$, when
$M^3$ is a closed connected orientable 3-manifold that admits a smooth
embedding into $\sharp_n\ncp$, the connected sum of complex projective planes
with negative orientation.  Ding's argument can be used to extend his result to
the cases where $\sharp_n\ncp$ is replaced by any definite closed smooth
4-manifold and also to include the 3-manifolds with boundary whose doubles have
the required embedding property.  Although the authors do not know of any
classification of 3-manifolds with respect to existence of the embeddings into
definite 4-manifolds, there are some well known examples of 3-manifolds, most
notably the Poincare sphere, that do not smoothly embed into a definite closed
4-manifold.

\proclaim{Theorem 4} Let $M^3$ be a compact 3-manifold.  Then  $M\times\r$ has
infinitely many different smooth structures.\endproclaim

\demo{Proof} We prove the theorem first for the case when $M$ is closed and 
orientable.  It is well known that any orientable closed compact 3-manifold
embeds into $\sharp_k(S^2\times S^2)$, for a $k$ large enough, see
[K]. Briefly, such a manifold bounds a 4-dimensional handlebody with a 0- and
2-handles and such that all attaching framings are even.  It is an easy
exercise in Kirby's link calculus to show that the double of this handlebody is
diffeomorphic to $\sharp_k(S^2\times S^2)$ (where $k$ is the number of the
2-handles in the handlebody).  So, we fix such $k$ for the given 3-manifold
$M$.  Obviously, a bicollared neighborhood of $M$ gives an embedding of
$M\times\r$ into $\sharp_k(S^2\times S^2)$.  From Theorem 1 it follows that
there is an integer $n_1$, such that \br{n_1} can not be smoothly embedded into
$\sharp_k(S^2\times S^2)$. We form an end-sum of $M\times\r$ with \br{n_1}.
Clearly $(M\times\r)\natural\br{n_1}$ can not be diffeomorphic with
$M\times\r$ since the former can not be smoothly embedded into
$\sharp_k(S^2\times S^2)$.  However, by its construction, \br{n_1} smoothly
embeds into $\sharp_{n_1}(S^2\times S^2)$.  Let $k_2=k+n_1$.  We fix an
embedding of $M\times\r$ into the first $k$ summands of
$\sharp_{k_2}(S^2\times S^2)$ and an embedding of \br{n_1} into the remaining
$n_1$ copies. We may assume that the both embedding miss the connecting tube
between the first $k$ and the last $n_1$ copies of $S^2\times S^2$.  Now, we
can perform the end-sum of $M\times\r$ and \br{n_1} ambiently inside
$\sharp_{k_2}(S^2\times S^2)$.  The next step is to select an integer $n_2$
large enough so that \br{n_2} does not embed into $\sharp_{k_2}(S^2\times
S^2)$. The new smooth structure,
$(M\times\r)\natural\br{n_1}\natural\br{n_2}=
(M\times\r)\natural\br{(n_1+n_2)}$ is different from the previous two and it
embeds into $\sharp_{k_3}(S^2\times S^2)$, where $k_3=k_2+n_2$.  We can
obviously iterate the construction to obtain a sequence of smooth structures,
$(M\times\r)\natural\br{n_i}$.

If $M$ is an orientable compact 3-manifold with the boundary we work with its
double, $DM$.  The double of $M$, and so $M$ itself, can be embedded into
$\sharp_k(S^2\times S^2)$ and the same construction produces a sequence of
smooth structures on the double.  We perform each end-sum on the same copy of
$M\times\r$ in the $DM\times\r$ and so we obtain a sequence of smooth
structures $(M\times\r)\natural\br{n_i}$ that are all standard near $\bd
M\times\r$.  Since each of these smooth manifolds requires a different minimal
number ``$n$'' to embed into $\sharp_n(S^2\times S^2)$, they are all different.

Finally, if $M$ is a non orientable compact 3-manifold it has a two-fold cover
$N$ that is an orientable 3-manifold.  $N$ smoothly embeds into
$\sharp_k(S^2\times S^2)$ for some $k$.  We choose $n_1$ as before, so that
\br{n_1} does not embed into $\sharp_k(S^2\times S^2)$.  We form the end-sum
$(M\times\r)\natural\br{n_1}$ and cover it by $(N\times\r)\natural_2\br{n_1}$.
By the construction $(N\times\r)\natural_2\br{n_1}$ does not embed into
$\sharp_k(S^2\times S^2)$ so it is not diffeomorphic to $N\times\r$.  Since a
diffeomorphism between $(M\times\r)\natural\br{n_1}$ and $M\times\r$ would lift
to their homotopically equivalent 2-fold covers it can not exist.  As in the 
previous
cases we obtain a sequence of different smooth structures on $M\times\r$.
\qed\enddemo

\remark{Remark} All the smooth structures on $M\times\r$ constructed above have
one of the two ends standard.  Obviously, by end-summing $M\times\r$ with
\br{n}'s on the both sides one can obtain doubly indexed countable family of
smooth structures.\endremark

\subhead{} Smooth structures on manifolds with collarable ends\endsubhead 
If an 4-manifold has countably many ends, each having a neighborhood of the
infinity that is {\it topologically} collarable as $S^3\times\r$, then Gompf
has shown that the manifold possesses uncountably many different smooth
structures [G2].  Ding [Di] has extended this results to include manifolds 
whose ends are topologically collared by $M^3\times\r$ where $M^3$
embeds into $\sharp_n\ncp$. 

\proclaim{Theorem 5} Let $X$ be an open topological four manifold with at 
least one topologically collarable end, i.e. an end homeomorphic to
$M\times\r$, for some closed 3-manifold $M$. If $X$ has more than one
end we assume that for one of the topologically collarable ends there
are only finitely many other ends that have this topological type.
Then, $X$ has infinitely many different smooth structures.
\endproclaim

Simple examples of manifolds to which the theorem applies are the interiors of
compact 4-manifolds with boundaries.

\demo{Proof} Let $X$ be a 4-manifold with a single end. Since any
open 4-manifold is smoothable we assume that $X$ has a smooth structure. Next
we form $X\natural \br{\infty}$.  We fix a homeomorphism $h$ between
$M\times\r$ and a neighborhood of the end of $X\natural \br{\infty}$ such that
the $n$th copy of $K_1$ is embedded into $M\times(n,n+1)$. For any $t\in\r, \
t\geq 1,$ we define $X_t$ as $(X\natural \br{\infty}) - h(M\times[t,+\infty))$.

We claim that if $t<s$ and $M\times(t,s)$ contains a copy of $K_1$ then $X_t$
and $X_s$ cannot be diffeomorphic.  From this it follows that $X$ has
infinitely many smooth structures.  To prove this claim we assume that there
are two positive real numbers, $t,\ s$, satisfying the properties above, but
such that $X_t$ is diffeomorphic to $X_s$. Then there is a neighborhood of the
end of $X_s$, denoted by $V$ that does not intersect the copy of
$K_1$ in $M\times(t,s)$.  The diffeomorphism between $X_s$ and $X_t$
identifies another copy of $V$ as a neighborhood of the end of $X_t$. We glue
to the end of $X_s$ infinitely many copies of the region between the two
copies of $V$ in $X_s$.  Using this smooth structure we find separating
smoothly embeddings of 3-manifolds: $M_1\hookrightarrow M\times(-1,0)$ and
$M_2 \hookrightarrow V\subset X_t$.  For $i=1,2$, $M_i$ bounds a smooth simply
connected spin 4-manifold $Z_i$.  We glue together $Z_1$, $Z_2$ and the piece
of $M\times(-1,t)$ bounded by $M_1\coprod M_2$. The resulting closed smooth
manifold $Z$ is spin (but not necessarily simply connected) and by adding
copies of the K3 surface with appropriate orientation we may assume that $Z$
has signature zero. Therefore $Z$ has the same intersection form as a
connected sum of hyperbolics.  Next we find $n$ large enough so that $K_n$
does not smoothly embed into a closed 4-manifold with the same intersection
form as $Z$. Note that we can embed $K_n$ in our newly constructed periodic
end by ambiently connecting together $n$ copies of $K_1$.  Also, by the
construction there is a copy $V$ passed the embedded $K_n$.  Now we can use
$M_2$ in this copy of $V$ to construct a smooth manifold homotopically
equivalent to $Z$ and containing $K_n$, which by our choice of $n$ is
impossible.

In the case that $X$ has more than one end we have assumed that at least one of
them is homeomorphic to $M\times\r$, for some 3-manifold $M$ and furthermore we
have assumed that there are finitely many ends of $X$ homeomorphic to
$M\times\r$. We repeat the construction from above on this end and obtain a
sequence of smooth structures $X_n$.  Since there can be only finitely many
diffeomorphisms that permute the ends homeomorphic to $M\times\r$ the sequence
$X_n$ contains infinitely many different smooth structures on $X$.\qed\enddemo

\subhead{}Final comments and questions\endsubhead 
The method used in the proof of Theorem 5 could also be used to prove Theorem 4
and part of Corollary 2.  However, the given proofs are more natural and
instructive.  Moreover, the proof of Theorem 4 can be used to show that
$M\times\r$ has infinitely many smooth structures if $M$ is an open 3-manifold
with all ends collarable.  Namely, one can trim all the ends to obtain a
compact 3-manifold with boundary.  Obviously Theorem 4 is true if $M$ is a
3-manifold such that $M\times\r$ is homeomorphic to $N\times\r$, where $N$ is a
3-manifold with collarable ends. An example of such an $M$ is the Whitehead
manifold whose end is not collarable, but whose product with the real line is
$\r^3\times\r$.  This leads to the question whether for any 3-manifold $M$ with
finitely generated homology there is a 3-manifold $N$ as above.

All our results produce only countable families of smooth structures.  However,
by shaving the end of any member of these families we can produce an
uncountable family of smooth structures.  Either there are uncountably many
different smooth types between them or the end of the original smooth structure
is periodic.  Although our methods fail to distinguish between them we
conjecture that the former is the case.

\enddocument